\documentclass[12pt]{amsart}
\usepackage{amssymb,latexsym,comment,url}

\newcommand{\EE}{{\mathcal E}}

\newcommand{\CC}{\operatorname{\mathcal C}}

\newcommand{\Q}{{\mathbb Q}}

\newcommand{\rank}{\operatorname{rank}}

\newcommand{\Z}{{\mathbb Z}}

\newcommand{\Magma}{{\sf MAGMA }}

\newenvironment{Proof}{\par\noindent{\sc Proof:}}%
                      {\hspace*{\fill}\nobreak$\Box$\par\medskip}
                       {\hspace*{\fill}\nobreak$\Box$\par\medskip}

\newtheorem{Proposition}{Proposition}[section]
\newtheorem{Theorem}[Proposition]{Theorem}

\newtheorem{Corollary}[Proposition]{Corollary}

\theoremstyle{definition}
\newtheorem{Definition}[Proposition]{Definition}

\renewcommand{\baselinestretch}{1.1}

\begin{document}

\title[Geometric Progressions on Hyperelliptic Curves]%
{On Geometric Progressions on Hyperelliptic Curves}

\author[M. Alaa]%
{Mohamed~Alaa}
\address{Department of Mathematics, Faculty of Science, Cairo University, Giza, Egypt}
\email{malaa@sci.cu.edu.eg}

\author[M. Sadek]%
{Mohammad~Sadek}
\address{American University in Cairo, Mathematics and Actuarial Science Department, AUC Avenue, New Cairo, Egypt}
\email{mmsadek@aucegypt.edu}

\begin{abstract}
Let $C$ be a hyperelliptic curve over $\Q$ described by $y^2=a_0x^n+a_1x^{n-1}+\ldots+a_n$, $a_i\in\Q$. The points $P_{i}=(x_{i},y_{i})\in C(\mathbb{Q})$, $i=1,2,...,k,$ are said to be in a geometric progression of length $k$ if the rational numbers $x_{i}$, $i=1,2,...,k,$ form a geometric progression sequence in $\Q$, i.e., $x_i=pt^{i}$ for some $p,t\in\Q$. In this paper we prove the existence of an infinite family of hyperelliptic curves on which there is a sequence of rational points in a geometric progression of length at least eight.
\end{abstract}

\maketitle


\section{Introduction}
Let $C$ be a hyperelliptic curve defined over the rational field $\Q$ by a hyperelliptic equation of the form $y^2=f(x)$, $\deg f(x)\ge3$. One may construct a sequence of $\Q$-rational points in $C(\Q)$ such that the $x$-coordinates of these rational points form a sequence of rational numbers which enjoys a certain arithmetic pattern. For instance, an arithmetic progression sequence on $C$ is a sequence $(x_i,y_i)\in C(\Q),\;i=1,2,\ldots,$ where $x_i=a+ib$ for some $a,b\in\Q$. In a similar fashion one may define a geometric progression sequence on $C$.

In \cite{Bremner} Bremner discussed arithmetic progression sequences on elliptic curves over $\Q$. He investigated the size of these sequences and he produced elliptic curves with long arithmetic progression sequences. His techniques were improved and used to generate infinitely many elliptic curves with long arithmetic progression sequences of rational points, see \cite{Campbell, Macleod, Ulas1}. In \cite{Ulas2} arithmetic progression sequences on genus 2 curves were considered.

  A certain family of algebraic curves was studied by Bremner and Ulas in \cite{BremnerUlas}. They proved the existence of an infinite family of algebraic curves defined by $y^2 = ax^n + b,\,n\ge 1,\, a, b \in \mathbb Q,$ with geometric progression sequences of rational points of length at least 4. They remarked that their method can be exploited in order to increase the length of these sequences to be 5.

  In this note we examine geometric progression sequences on hyperelliptic curves. We start with proving that unlike geometric progressions on the rational line, geometric progression sequences on hyperelliptic curves are finite. A certain family of hyperelliptic curves defined by an equation of the form  $y^{2}=ax^{2n}+bx^{n}+a,\;n \in \mathbb{N},\; a,b \in \mathbb{Q},$ is displayed. Each hyperelliptic curve in this family possesses a geometric progression sequence of rational points whose length is at least 8. In fact, those hyperelliptic curves are parametrized by an elliptic surface $\mathcal H$ with positive rank. In particular, to each point of infinite order on $\mathcal H$ one can associates a hyperelliptic curve with a geometric progression sequence of length at least $8$.

  It is worth mentioning that other types of sequences of rational points on algebraic curves are being studied. For example, in \cite{SadekKamel} an infinite family of elliptic curves is constructed such that every elliptic curve in the family has a sequence of rational points whose $x$-coordinates form a sequence of consecutive rational squares. The length of the latter sequence is at least $5$.

\section{Geometric progression sequences on hyperelliptic curves}
Let $C$ be a hyperelliptic curve defined over a number field $K$ by the equation $y^2=f(x)$ where $f(x)\in K[x]$ is of degree $n\ge 3$, and $f(x)$ has no double roots. The set $C(K)$ of $K$-rational points on $C$ is defined by $\displaystyle C(K)=\{(x,y):y^2=f(x),\,x,y\in K\}$.
\begin{Definition}
Let $C:y^2=f(x)$ be a hyperelliptic curve over a number field $K$. The sequence $P_i=(x_i,y_i)\in C(K),\,i=1,2,\ldots,$ is said to be a {\em geometric progression sequence} in $C(K)$ if there are $p,t\in K^{\times}$ such that $x_i=pt^i$. In other words, the $x$-coordinates of the rational points $P_i$ form a geometric progression sequence in $K$.
\end{Definition}

We assume throughout that our geometric progression sequences contain distinct rational points, in particular $t\not\in\{\pm1\}$.

We will show that unlike geometric progressions in $K$, geometric progression sequences in $C(K)$ are finite.

\begin{Theorem}
\label{thm:length}
Let $C:y^2=f(x)$ be a hyperelliptic curve over a number field $K$ with $\deg f(x)\ge 3$. Let $(x_i,y_i)$ form a geometric progression sequence in $C(K)$. Then the sequence $(x_i,y_i)$ is finite.
\end{Theorem}
\begin{Proof}
If $\deg f(x)\ge 5$, then it follows that the genus $g$ of $C$ satisfies $g\ge 2$. In view of Faltings' Theorem, \cite{Falting}, $C(K)$ is finite.

If $\deg f(x)=3$ or $4$, then $C$ is an elliptic curve. Assume that $f(x)=a_0x^4+a_1x^3+a_2x^2+a_3x+a_4$. Assume on the contrary that there is an infinite sequence $(x_i,y_i)\in C(K)$, $x_i=pt^i,\,i=1,2,\ldots,$ for some $p,t\in K^{\times}$. Considering the subsequence $(x_{2i},y_{2i}),$ $i=1,2,\ldots$, one obtains
\[y_{2i}^2=a_0p^4t^{8i}+a_1p^3t^{6i}+a_2p^2t^{4i}+a_3p t^{2i}+a_4,\,i=1,2,\ldots.\]
In particular, the rational points $(t^i,y_i),i=1,2,\ldots,$ form an infinite sequence of rational points on the new hyperelliptic curve
\[C':y^2=a_0p^4x^8+a_1p^3x^6+a_2p^2x^4+a_3px^2+a_4.\]
This contradicts Faltings' Theorem, since the genus of $C'$ is $2$ if $a_0=0$; $3$ if $a_0\ne 0$.
\end{Proof}

The theorem above motivates the following definition. Given a geometric progression sequence $(x_i,y_i)$, $i=1,2,\ldots,k,$ in $C(K)$, the positive integer $k$ will be called the {\em length} of the sequence.

\section{Hyperelliptic curves with long geometric progressions}

In this note, we consider the family of hyperelliptic curves over $\Q$ described by the equation $y^2=ax^{2n}+bx^n+a$ where $a,b\in\Q$, and $n\ge 2$. We introduce an infinite family of these hyperelliptic curves with geometric progression sequences of length at least $8$. We remark that the existence of a sequence of rational points $(t^i,y_i)$, $i=1,2,\ldots,k$, in geometric progression on one of these hyperelliptic curves is equivalent to the existence of the following geometric progression sequence of rational points $(t^{ni},y_i)$ on the conic $y^2=ax^2+bx+a$. In fact, we will establish the existence of such an infinite family of conics on which there exist geometric progression sequences of rational points whose $x$-coordinates are $t^{-7},t^{-5},t^{-3},t^{-1},t,t^3,t^5,t^7$ for some $t\in\Q\setminus\{-1,0,1\}$.

We start with assuming that the points $(t,U)$ and $(t^3,V)$ are two rational points in $C(\Q)$ where $C$ is given by $y^2=f(x)=ax^2+bx+a$. This implies that
\begin{eqnarray*}
U^{2}&=&at^{2}+bt+a,\\
V^{2}&=&at^{6}+bt^{3}+a,
\end{eqnarray*}
hence
\begin{eqnarray}
\label{eq1}
a&=&\frac{t^{2}U^{2}-V^{2}}{(t^{2}-1)^{2}(t^{2}+1)},\nonumber\\
b&=&\frac{(t^{4}-t^{2}+1)U^{2}-V^{2}}{t(t^{2}-1)^{2}}.
\end{eqnarray}

From the symmetry of the polynomial $f(x)$, one observes that if the points $(t,U)$ and $(t^{3},V)$ are in $C(\Q)$, then so are the points $\displaystyle(t^{-1},Ut^{-1})$ and $(t^{-3},Vt^{-3})$. So we already have four points in geometric progression in $C(\Q)$.

In order to increase the length of the progression, we assume that $(t^{5},R)$ is in $C(\Q)$, hence $(t^{-5},Rt^{-5})$ is in $C(\Q)$ as well. Given the description of $a$ and $b$, (\ref{eq1}), one obtains
$$R^{2}=-t^{2}(t^{4}+1)U^{2}+(1+t^{2}+t^{4})V^{2}.$$
\begin{Theorem}
\label{Thm2}
The conic $\mathcal{C}:R^{2}=-t^{2}(t^{4}+1)U^{2}+(1+t^{2}+t^{4})V^{2}$ defined over $\Q(t)$ has infinitely many rational points given by the following parametrization
\begin{eqnarray}
\label{eq2}
U&=&t^{2}(1+t^{4})p^{2}+(1+t^{2}+t^{4})q^{2}-2t(1+t^{2}+t^{4})pq,\nonumber\\
V&=&t^2(1+t^4)p^2+(1+t^{2}+t^{4})q^{2}-2t(1+t^{4})pq,\nonumber\\
R&=&t^{3}(1+t^{4})p^{2}-t(1+t^{2}+t^{4})q^{2}.
\end{eqnarray}
\end{Theorem}
\begin{Proof}
The point $(U:V:R)=(1:t:t^2)$ lies in $\mathcal{C}(\Q(t))$. This implies the existence of infinitely many rational points on the conic $\CC$. These rational points are given by a parametric description, and the parametrization can be found in \cite[p. 69]{Mordell}.
\end{Proof}
\begin{Corollary}
There exists an infinite family of conics $y^2=ax^2+bx+a$, $a,b\in\Q$, containing 6 rational points in geometric progression. In particular, there exist infinitely many hyperelliptic curves described by the equation $y^2=ax^{2n}+bx^n+a$ with 6 rational points in geometric progression.
\end{Corollary}

In what follows we parametrize the family of conics $C:y^2=ax^2+bx+a$ containing a seventh rational point $(t^7,S)$. We recall that the existence of this seventh rational point implies the existence of an eighth point $(t^{-7},St^{-7})$ on the conic $C$. The point $(t^7,S)$ satisfies the equation of the conic where $a,b$ are described as in (\ref{eq1}) and (\ref{eq2}). This gives rise to the following curve over $\Q(t)$
\begin{multline}
\label{eq3}
\mathcal{H}:S^2=H_t(p,q):=t^8(1+t^4)^2p^4+4t^5(1+2t^4+2t^8+t^{12})p^3q-2t^4(4+3t^2+9t^4+4t^6+9t^8+3t^{10}+4t^{12})p^2q^2\\
        +4t^3(1-t^2+t^4)(1+t^2+t^4)^2pq^3+t^4(1+t^2+t^4)^2q^4.
        \end{multline}
        \begin{Theorem}
        \label{thm3}
        The curve $\mathcal H$ defined over $\Q(t)$ is an elliptic curve with $\rank \mathcal{H}(\Q(t))\ge 1$.
        \end{Theorem}
        \begin{Proof}
        The following point lies in $\mathcal H(\Q(t))$ {\footnotesize$$(p:q:S)=\left(\frac{-t}{1-t^2+t^4}:1-\frac{3+2t^2+4t^4+2t^6+3t^8}{2(1+t^4+t^8)}:\frac{t^2(3+4t^2+8t^4+8t^6+10t^8+8t^{10}+8t^{12}+4t^{14}+3 t^{16})}{4(1-t^2+t^4)^2(1+t^2+t^4)}\right).$$}
The existence of the latter rational point in $\mathcal H(\Q(t))$ implies that the curve $\mathcal H$ is birationally equivalent over $\Q(t)$ to its Jacobain $\mathcal E$ described by $Y^2=4X^3-g_2X-g_3$ where
\begin{eqnarray*}
g_2&=&\frac{4}{3}t^8(1+t^2+t^4)^2(1+t^2+4t^4+t^6+7t^8+t^{10}+4t^{12}+t^{14}+t^{16}),\\
g_3&=&-\frac{4}{27}t^{12}(1+t^2+t^4)^4(2+t^2+3t^4+15t^6-9t^8+30t^{10}-9t^{12}+15t^{14}+3t^{16}+t^{18}+2t^{20}),
\end{eqnarray*}
see \cite{Mordell}. The point $P=(X_P,Y_P)$ where
{\footnotesize\begin{eqnarray*}X_P&=&-\frac{t^4(1+t^2+t^4)^2(2-5t^2-2t^4-2t^6-2t^8-5t^{10}+2t^{12})}{3(1+t^2)^4},\\
Y_P&=&\frac{4t^7(1+t^2+t^4)^2}{(1+t^2)^6}(1+t^2+2t^4+2t^6+3t^8+2t^{10}+3t^{12}+2t^{14}+2t^{16}+t^{18}+t^{20})\end{eqnarray*}}
is a point in $\mathcal E(\Q(t))$. In fact, specializing $t=2$ and using \Magma, \cite{MAGMA}, we find that the specialization of the point $P$ is a point of infinite order on the specialization of $\EE$ when $t=2$. It follows that the point $P$ itself is a point of infinite order in $\EE(\Q(t))$.
\end{Proof}
\begin{Corollary}
\label{cor:hyperelliptic}
Fix $t_0\in\Q$. For any nontrivial geometric progression sequence of the form $t_0^{\pm1},t_0^{\pm3},t_0^{\pm5},t_0^{\pm 7}$, there exist infinitely many hyperelliptic curves $C_m:y^2=a_mx^{2n}+b_mx^n+a_m,\;m\in\Z\setminus\{0\},n\ge2,$ such that the numbers $t_0^{\pm i},i=1,3,5,7,$ are the $x$-coordinates of rational points on $C_m$.
\end{Corollary}
\begin{Proof}
The point $P=(p:q:S)$ described by {\footnotesize$$\left(\frac{-t_0}{1-t_0^2+t_0^4}:1-\frac{3+2t_0^2+4t_0^4+2t_0^6+3t_0^8}{2(1+t_0^4+t_0^8)}:\frac{t_0^2(3+4t_0^2+8t_0^4+8t_0^6+10t_0^8+8t_0^{10}+8t_0^{12}+4t_0^{14}+3 t_0^{16})}{4(1-t_0^2+t_0^4)^2(1+t_0^2+t_0^4)}\right)$$} is a point of infinite order on the curve $\mathcal H$ over $\Q(t_0)$, see Theorem \ref{thm3}. For any nonzero $m$, we write $mP=(p_m:q_m:S_m)$ for the $m$-th multiple of $P$.

Substituting these values of $p_m,q_m\in\Q(t_0)$ into (\ref{eq2}), one obtains a parametric solution $U_m,V_m,R_m$ of the quadratic $R^2=-t_0^2(t_0^4+1)U^2+(1+t_0^2+t_0^4)V^4$. Hence, one obtains $a_m$ and $b_m$ by substituting $U_m$ and $V_m$ into the formulas of $a,b$ in (\ref{eq1}).

   We get an infinite family of hyperelliptic curves $C_m:y^2=a_mx^{2n}+b_mx^n+a_m$, where $m$ is nonzero. This family satisfies the property that the points $(t_0^i,u_i),(t_0^{-i},u_it_0^{-i})$, $i=1,3,5,7$, are lying in $C_m(\Q)$ for some $u_i\in\Q$. Thus, one obtains an infinite family of hyperelliptic curves with an $8$-term geometric progression sequence of rational points.
  \end{Proof}
\section{A numerical example}
  The curve $C: y^{2}=a(T)x^{2n}+b(T)x^{n}+a(T),n \in \mathbb N$, where  $a(T)$ is given by \\ 
{\footnotesize  $\displaystyle \frac{T^{4n}(1+T^{2n})(1+T^{8n})}{2(1+T^{4n})(-1+T^{2n}-T^{4n}+T^{6n}-T^{8n}+T^{10n})^2}$} and $b(T)$ is defined by {\footnotesize$$\frac{1-2T^{2n}-T^{4n}-12T^{6n}-3T^{8n}-14T^{10n}-13T^{12n}-40T^{14n}-13T^{16n}-14T^{18n}-3T^{20n}-12T^{22n}-T^{24n}-2T^{26n}+T^{28n}}{16T^{3n}
(-1+T^{2n})^2(1+T^{4n})^2(1-T^{2n}+T^{4n})^2(1+T^{2n}+T^{4n})^2},$$} has the following $8$-term geometric progression sequence

{\footnotesize \begin{gather*}
\left(T^{-7},\frac{3+4T^{2n}+5T^{4n}+4T^{6n}+5T^{8n}+4T^{10n}+3T^{12n}}{4T^{5n}(1+2T^{4n}+2T^{8n}+T^{12n})}\right),\\
\left(T^{-5},\frac{1+4T^{2n}+3T^{4n}+4T^{6n}+3T^{8n}+4T^{10n}+T^{12n}}{4T^{4n}(1+2T^{4n}+2T^{8n}+T^{12n})}\right),
\left(T^{-3},\frac{1+3T^{4n}+4T^{6n}+3T^{8n}+T^{12n}}{4T^{3n}(1+2T^{4n}+2T^{8n}+T^{12n})}\right),\\
\left(T^{-1},\frac{-1+T^{4n}+4T^{6n}+T^{8n}-T^{12n}}{4T^{2n}(1+2T^{4n}+2T^{8n}+T^{12n})}\right),
\left(T,\frac{-1+T^{4n}+4T^{6n}+T^{8n}-T^{12n}}{4T^n(1+2T^{4n}+2T^{8n}+T^{12n})}\right),\\
\left(T^3,\frac{1+3T^{4n}+4T^{6n}+3T^{8n}+T^{12n}}{4(1+2T^{4n}+2T^{8n}+T^{12n})}\right),
\left(T^5,\frac{T^n(1+4T^{2n}+3T^{4n}+4T^{6n}+3T^{8n}+4T^{10n}+T^{12n})}{4(1+2T^{4n}+2T^{8n}+T^{12n})}\right),\\
\left(T^7,\frac{T^{2n}(3+4T^{2n}+5T^{4n}+4T^{6n}+5T^{8n}+4T^{10n}+3T^{12n})}{4(1+2T^{4n}+2T^{8n}+T^{12n})}\right).
\end{gather*}}

For example, When $n=2$ and $t=2$, one has the elliptic curve
 \[y^2=\frac{142608512}{250308167443425}x^4+\frac{62553486161362657 }{65873099809751270400}x^2+\frac{142608512}{250308167443425}\] which contains the following $8$-term geometric progression sequence
$$(2^{-7},\frac{54871363}{69258448896}),(2^{-5},\frac{21185345}{17314612224}),(2^{-3},\frac{5663659}{1442884352}),(2^{-1},\frac{16695041}{1082163264}),$$
$$(2,\frac{16695041}{270540816}),(2^3,\frac{5663659}{22545068}),(2^5,\frac{21185345}{16908801}),(2^7,\frac{219485452}{16908801}).$$

\section{A remark on geometric progressions of length 10}
In order to extend the length of the $8$-term geometric progression sequence we constructed in Corollary \ref{cor:hyperelliptic} to a geometric progression of length $10$, one assumes that a point of the form $(t^9,S')$, and consequently the point $(t^{-9},S't^{-9})$, exists on the hyperelliptic curve $y^2=a(t)x^{2n}+b(t)x^n+a(t)$. This yields the existence of a rational point $(p:q:S')$ on the elliptic curve $\mathcal L$ defined by
{\footnotesize\begin{multline*}
S'^2=H_t'(p,q):=t^{10}(1+t^4)^2p^4+4t^5(1+t^2+2t^4+2t^6+2t^8+2t^{10}+2t^{12}+t^{14}+t^{16})p^3q\\
-2t^4(4+6t^2+11t^4+11t^6+12t^8+11t^{10}+11t^{12}+6t^{14}+4t^{16})p^2q^2\\
+4t^3(1+2t^2+3t^4+3t^6+3t^8+3t^{10}+3t^{12}+2t^{14}+t^{16})pq^3
+t^6(1+t^2+t^4)^2q^4.
\end{multline*}}
One recalls that the pair $(p,q)$ makes up the first two coordinates of a point $(p:q:S)$ on the elliptic curve $\mathcal H:S^2=H_t(p,q)$ defined over $\Q(t)$. This implies that one needs to find a solution $(p,S,S')$ on the genus $5$ curve $\mathcal C$ defined by the affine equation
\[S^2=H_t(p,1),\;S'^2=H'_t(p,1).\] In view of Faltings' Theorem, a genus five curve possesses finitely many rational points. Therefore, one reaches the following result.
\begin{Proposition}
 Fix $t_0\in\Q$. For any nontrivial $10$-term geometric progression sequence of the form $t_0^{\pm1},t_0^{\pm3},t_0^{\pm5},t_0^{\pm 7},t_0^{\pm 9}$, there exist finitely many hyperelliptic curves of the form $C:y^2=ax^{2n}+bx^n+a,\,a,b\in\Q,$ such that the numbers $t_0^{\pm i},i=1,3,5,7,9,$ are the $x$-coordinates of rational points in $C(\Q).$
\end{Proposition}

\hskip-12pt\emph{\bf{Acknowledgements.}}
We would like to thank Professor Nabil Youssef, Cairo University, for his support, careful reading of an earlier draft of the paper, and several useful suggestions that improved the manuscript.

\end{document}